\documentclass[12pt]{amsart}
\usepackage{amsmath}
\usepackage{amsfonts, amssymb, euscript, mathrsfs}
\usepackage{amsthm, upref}
\usepackage{graphicx}
\usepackage[usenames, dvipsnames]{color}
\usepackage{float}
\usepackage{pdfpages}
\usepackage[caption=false]{subfig}

\usepackage{enumerate}
\usepackage[margin=1in]{geometry}
\usepackage{aliascnt}
\usepackage{hyperref}
\usepackage{verbatim}

\allowdisplaybreaks[4]

\newcommand{\BR}{\mathbb{R}}

\newcommand{\SL}{\sum\limits}

\newcommand{\de}{\delta}

\newcommand{\CA}{\mathcal A}
\newcommand{\CL}{\mathcal L}

\newcommand{\CN}{\mathcal N}

\newcommand{\si}{\sigma}

\newcommand{\pa}{\partial}
\renewcommand{\phi}{\varphi}

\newcommand{\la}{\lambda}

\newcommand{\ol}{\overline}

\newcommand{\norm}[1]{\lVert#1\rVert}
\renewcommand{\comment}[1]{}

\newcommand{\mP}{\mathbf{p}}

\newcommand{\md}{\mathrm{d}}

\DeclareMathOperator{\Exp}{Exp}

\DeclareMathOperator{\mes}{mes}


\begin{document}


\theoremstyle{plain}
\newtheorem{thm}{Theorem}[section]
\newtheorem{lemma}[thm]{Lemma}
\newtheorem{prop}[thm]{Proposition}
\newtheorem{cor}[thm]{Corollary}
\newtheorem{open}[thm]{Open Problem}

\theoremstyle{definition}
\newtheorem{defn}{Definition}
\newtheorem{asmp}{Assumption}
\newtheorem{notn}{Notation}
\newtheorem{prb}{Problem}

\theoremstyle{remark}
\newtheorem{rmk}{Remark}
\newtheorem{exm}{Example}
\newtheorem{clm}{Claim}

\author{Clayton Barnes}

\address{\scriptsize Technion--Israel Institute of Technology}

\email{clayleroy2@gmail.com}

\author{Andrey Sarantsev}

\address{\scriptsize{Department of Mathematics and Statistics, University of Nevada, Reno}} 

\email{asarantsev@unr.edu}

\title[A Note on Jump Atlas Models]{A note on Jump Atlas Models}  

\date{\today}



\begin{abstract} The market weight of a stock is its capitalization (cap) divided by the total market cap.  Rank these weights from top to bottom. The capital distribution curve is a plot of weights versus ranks. For the US stock market, it is linear on a double logarithmic scale, and stable with respect to time (Fernholz, 2002). This property has been captured by models with rank-dependent dynamics: Each stock's cap logarithm is a Brownian motion with drift and diffusion coefficients depending on its current rank (Chatterjee, Pal, 2010). However, short-term stock movements have heavy tails. One can add jumps to Brownian motions to capture this. Observed time stability follows from a long-term stability result, stated and proved here. Via simulations, we find which properties of continuous models are preserved after adding jumps. 
\end{abstract}

\maketitle

\thispagestyle{empty}

\section{Introduction}

\subsection{Motivation and definitions} A {\it market model} is a collection $(S_1, \ldots, S_N)$ of $N$ positive real-valued continuous-time random processes $S_i = (S_i(t),\, t \ge 0)$: The {\it capitalization} of the $i$th stock at time $t \ge 0$ is given by $S_i(t)$. The total capitalization of the stock market is $S(t) = S_1(t) + \ldots + S_N(t)$, and the {\it market weight} of the $i$th stock at time $t \ge 0$ is
\begin{equation}
\label{eq:weights}
\mu_i(t) = \frac{S_i(t)}{S(t)}.
\end{equation}
Rank these market weights from top to bottom: $\mu_{(1)}(t) \ge \ldots \ge \mu_{(N)}(t)$. \cite[Chapter 5, page 95]{FernholzBook} contains the double logarithmic plot of ranked market weights $(\ln n, \ln\mu_{(n)}(t)),\, n = 1, \ldots, N$, December 31, 1929--1999 (8 plots, every 10 years), stocks from the Center of Research in Securities Prices at the University of Chicago. This plot exhibits two features:

\begin{enumerate}[(A)]
\item It is stable over time: It is almost the same for all $t$.
\item It is close to a straight line, except at lower and upper ends.
\end{enumerate}

\begin{figure}[t]
\includegraphics[width = 7cm, height = 5cm]{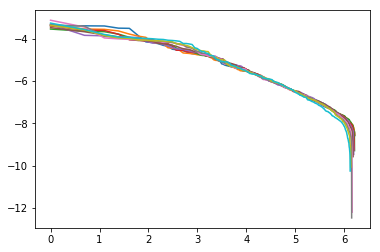}
\caption{Double logarithmic plot $(\ln k,\, \ln\mu_{(k)}(t)),\, k = 1, \ldots, N$, of real-world ranked market weights, Dec 31, 2009--2018}
\label{fig:SPplot}
\end{figure}

We reproduce this (and see the same properties) in Figure~\ref{fig:SPplot} for the S\&P 500 stocks included in this index as of May 9, 2019, in December 31, 2009--2018, using the data from the \texttt{YCharts} database. (The code is given on \texttt{GitHub}.) These observations (a) and (b) were explained by the following rank-based model in \cite{BFK2005}. 

We introduce a few pieces of notation. The Dirac delta measure at point $b$ is denoted by $\delta_b$. For a vector $x = (x_1, \ldots, x_N) \in \mathbb R^N$, its {\it ranking permutation} $\mathbf{p}_x$ for the vector $x$ is defined as the unique permutation on $\{1, \ldots, N\}$ such that:
\begin{itemize}
\item $x_{\mathbf{p}_x(i)} \le x_{\mathbf{p}_x(j)}$ for $i < j$ (this permutation is ranking the vector $x$).
\item if $i < j$ and $x_{\mathbf{p}_x(i)} = x_{\mathbf{p}_x(j)}$, then $\mathbf{p}_x(i) < \mathbf{p}_x(j)$ (ties resolved in lexicographic order).
\end{itemize}
We rank from bottom to top, except market weights, ranked from top to bottom. 

\begin{defn} Take $N$ continuous adapted random processes $X_i = (X_i(t),\, t \ge 0),\, i = 1, \ldots, N$. Rank them at each time $t$ from bottom to top: $X_{(k)}(t) = X_{\mathbf{p}_{X(t)}(k)}$ for $k = 1, \ldots, N$; $X_{(1)}(t) \le \ldots \le X_{(N)}(t)$. If they satisfy a system of stochastic differential equations:
\begin{equation}
\label{eq:CBP}
\mathrm{d}X_i(t) = \sum\limits_{k=1}^N1(X_i(t)\ \mbox{has rank}\ k\ \mbox{at time}\ t)\,(g_k\,\mathrm{d}t  + \sigma_k\,\mathrm{d}W_i(t)),\ i = 1, \ldots, N,
\end{equation}
with $W_1, \ldots, W_N$ independent Brownian motions, and constant numbers $g_k$ and $\sigma_k$, then $(X_1, \ldots, X_N)$ is called a {\it system of competing Brownian particles;} with {\it drift and diffusion coefficients for rank $k$} equal to $g_k$ and $\sigma_k^2$. The {\it gap process} is defined as
\begin{equation}
\label{eq:gap}
Z(t) = \left(X_{(2)} - X_{(1)}(t), \ldots, X_{(N)} - X_{(N-1)}(t)\right),\ t \ge 0. 
\end{equation}
\end{defn}
Weak existence and uniqueness in law for this system is proved in \cite{BFK2005}. Take such a system of competing Brownian particles, and consider the following market model:
\begin{equation}
\label{eq:basic}
S_i(t) = \exp(-X_i(t)),\quad i = 1, \ldots, N;\ t \ge 0.
\end{equation}

\begin{defn} 
A {\it stationary distribution,} or {\it invariant measure,} for the gap process $Z$, is any probability measure $\Pi$ in the orthant $[0, \infty)^{N-1}$ such that $Z(0) \sim \Pi$ implies $Z(t) \sim \Pi$ for every $t \ge 0$. This measure $\Pi$ is sometimes called a {\it stationary gap distribution.} The gap process is called {\it stable} if it has a unique stationary distribution $\Pi$, and we have: 
\begin{equation}
\label{eq:ergodic}
\sup\limits_{E \subseteq [0, \infty)^{N-1}}\left|\mathbb P(Z(t) \in E\mid Z(0) = z) - \Pi(E)\right| \to 0\ \mbox{as}\ t \to \infty.
\end{equation}
\end{defn}

We introduce the following notation. For any vector $x = (x_1, \ldots, x_N) \in \mathbb R^N$, we define two other vectors: {\it running average} $\tilde{x} \in \mathbb R^N$, and {\it centered vector} $\overline{x} \in \mathbb R^N$:
\begin{equation}
\label{eq:centering-procedure}
\tilde{x}_k := \frac1k(x_1 + \ldots + x_k);\quad \overline{x}_k := x_k - \tilde{x}_N,\quad k = 1, \ldots, N.
\end{equation}
From~\eqref{eq:centering-procedure}, we have the following algebraic relation:
\begin{equation}
\label{eq:relation}
\overline{x}_1 + \ldots + \overline{x}_k = k(\tilde{x}_k - \tilde{x}_N),\quad k = 1, \ldots, N.
\end{equation}
Assume the average drift $\tilde{g}_k$ for the bottom-ranked $k$ particles is greater than the average drift $\tilde{g}_N$ for all $N$ particles.
\begin{equation}
\label{eq:stable-CBP}
\tilde{g}_k > \tilde{g}_N, \quad k = 1, \ldots, N-1.
\end{equation}
Under this condition~\eqref{eq:stable-CBP}, it was proved in \cite{BFK2005, 5people} that the gap process is stable. This has the following intuition: The drift of the ``cloud'' composed of $k$ bottom-ranked particles is larger than the drift for all particles; therefore, the bottom-ranked $k$ particles and the top-ranked $N-k$ particles cannot eventually separate into two ``clouds''. Multiplying the numerator and the denominator of~\eqref{eq:weights} by $\exp(X_{(1)}(t))$, we see that the ranked market weights from~\eqref{eq:weights} can be computed from this gap process in~\eqref{eq:gap}, and vice versa. Thus the vector of ranked market weights
\begin{equation}
\label{eq:ranked-market}
\left(\mu_{(1)}, \ldots, \mu_{(N)}\right)
\end{equation}
has a unique stationary distribution $\mathcal{M}$ (which is a push-forward of $\mathcal P$), and a stability result similar to~\eqref{eq:ergodic} holds. For the case $\sigma_k^2 = \alpha k + \beta$ (for some $\alpha, \beta$), the distribution $\Pi$ is a product of exponentials, \cite{5people}. The exponential distribution with rate $\la$ is denoted by $\Exp(\la)$; and for $\sigma_k = 1,\, k = 1, \ldots, N$, this distribution takes the form
$$
\Pi = \bigotimes\limits_{k=1}^{N-1}\Exp\left(2k(\overline{g}_k - \overline{g}_N)\right).
$$
For the {\it Atlas model} (named after the ancient Atlas hero, holding the sky on his shoulders) 
\begin{equation}
\label{eq:Atlas}
g_k = 
\begin{cases}
g,\ k = 1;\\
0,\ k = 2, \ldots, N;
\end{cases}
\quad 
\sigma_k = 1,\ k = 1, \ldots, N;\quad g > 0,
\end{equation}
the limiting shape (as $N \to \infty$) of the double logarithmic plot of ranked market weights $(\ln k,\, \ln\mu_{(k)}(t)),\, k = 1, \ldots, N$ from~\eqref{eq:ranked-market} is essentially a straight line, see \cite{CP2010}. (The precise statement is a bit more complicated.) 

However, this model~\eqref{eq:CBP} has one disadvantage: Fluctuations of real-world stock prices are not well described by Gaussian distributions. Instead, these prices make occasional large jumpsm and fluctuations have heavy tails. This is captured by adding jump components to random processes describing the stock market dynamics. An important class of stock price models contains L\'evy processes, which are generalizations of a Brownian motion: $W_i(t) - W_i(s)$ for $s < t$ can have some distribution other than Gaussian. In this case, the trajectories are not continuous. (A Brownian motion is the only example of a L\'evy process with continuous trajectories.) Here is an example of a L\'evy process.

\begin{defn} A  {\it compound Poisson process} $X(t) = V_1 + \ldots + V_{N(t)}$, where $V_1, V_2, \ldots$ are independent identically distributed random variables, and $N(\cdot)$ is a {\it Poisson process} with {\it rate} $\lambda > 0$, independent of $V_1, V_2, \ldots$: That is, $N(t) - N(s)$ for $t > s$ is distributed as a Poisson random variable with mean $\lambda(t-s)$, and is independent of $N(u),\, 0 \le u \le s$. If $\mathcal Q$ is the distribution of $Z_i$, then $\Lambda := \lambda\mathcal Q$ is called the {\it spectral measure}. 
\end{defn}

The sum of a compound Poisson process and a Brownian motion with drift $g$ and diffusion $\sigma^2$ is also a L\'evy process. In fact, every L\'evy process with finitely many jumps on a finite time interval can be represented in this way. We say that this L\'evy process is {\it associated with triple} $(g, \sigma, \Lambda)$. 

\begin{defn} For every $k = 1, \ldots, N$, let $L_k$ be a L\'evy process associated with the triple $(g_k, \sigma_k, \Lambda_k)$, $k = 1, \ldots, N$; that is, we can express each $L_k$ as follows:
\begin{equation}
\label{eq:explicit-expression}
L_k(t) = L_k(0) + g_kt + \sigma_kW_k(t) + \sum\limits_{j=1}^{N_k(t)}V_j,\ t \ge 0,
\end{equation}
where $W_k$ is a Brownian motion, $N_k$ is a Poisson process with rate $\lambda_k := \Lambda_k(\mathbb R)$, independent of $W_k$, and $V_j \sim \lambda_k^{-1}\Lambda_k,\, j = 1, 2, \ldots$ are i.i.d. random variables, independent of $W_k$ and $N_k$. The processes $L_1, \ldots, L_N$ are assumed to be independent. Modifying~\eqref{eq:CBP} as follows:
\begin{equation}
\label{eq:CLP}
\mathrm{d}X_i(t) = \sum\limits_{k=1}^N1(X_i(t)\ \mbox{has rank}\ k\ \mbox{at time}\ t)\,\mathrm{d}L_k(t),\ i = 1, \ldots, N,
\end{equation}
we get a {\it system of competing L\'evy particles}. 
\end{defn}

In this article, we have the following goals:

\begin{enumerate}
\item To prove a long-term stability result for the gap process~\eqref{eq:gap} of~\eqref{eq:CLP};
\item To find an explicit stationary distribution $\Pi$ for the gap process, and the corresponding stationary distribution $\mathcal R$ for the vector of ranked market weights~\eqref{eq:ranked-market};
\item To show that for some choice of parameters $(g_k, \sigma_k, \Lambda_k)$, the capital distribution curve for a large number $N$ of companies is close to the straight line.
\end{enumerate}

We accomplish Goal 1 by proving a (slightly more general) result: Theorem~\ref{thm:first}. We accomplish Goal 3 by simulations of these particle systems. We are unable to accomplish Goal 2 in exact form: It was shown in \cite{Bad} that a stationary distribution for a reflected Brownian motion with jumps in the orthant (of which the gap process is a particular case) does not have a product form. Instead, we perform numerical simulations, find empirical properties of the stationary gap distribution, and compare it to models without jumps. 

Often, jumps in financial models are modeled using L\'evy processes with infinitely many jumps on a finite time interval, for example 
{\it $\alpha$-stable processes,} see \cite{ContBook}. Such L\'evy process cannot be represented as a sum of a Brownian motion and a compound Poisson process. It is important to establish results for these models. But the fluctuations of such processes do not have a second moment, and this makes modifying the proof in Section 3 quite difficult. 

\subsection{The main stability result} The condition~\eqref{eq:stable-CBP} for a system with jumps takes the form:
\begin{align}
\label{eq:stable-CLP}
\begin{split}
\tilde{m}_k &> \tilde{m}_N,\ k = 1, \ldots, N-1,\\
m_k := g_k + f_k, & \quad f_k := \int_{\mathbb R}z\Lambda_k(\mathrm{d} z),\quad k = 1, \ldots, N,
\end{split}
\end{align}
and $\tilde{m}_k$ is defined from $m_k$ using~\eqref{eq:centering-procedure}. Using~\eqref{eq:relation}, we can rewrite the condition~\eqref{eq:stable-CLP} as
\begin{equation}
\label{eq:stability-new}
\min\limits_{k = 1, \ldots, N-1}(\overline{m}_1 + \ldots + \overline{m}_k) > 0,\quad \mbox{where}\quad \overline{m}_k = \overline{g}_k + \overline{f}_k,\ \ k = 1, \ldots, N.
\end{equation}
Each $m_k$ has the following meaning: For every $t \ge 0$, $\mathbb E\left[L_k(t) - L_k(0)\right] = m_kt$. This {\it effective drift} $m_k$ contains two terms: the actual drift $g_k$ in~\eqref{eq:explicit-expression}, and the product $f_k$ of $\lambda_k$, the rate of the Poisson process $N_k$, with the expected value of $Z_j$, distributed as $\lambda_k^{-1}\Lambda_k$. Therefore, $\tilde{m}_k$ is the average drift for the bottom-ranked $k$ particles out of $X_1, \ldots, X_N$: The drift for the ``cloud'' of these particles, as they move together. Thus, $\tilde{m}_k$ plays the same role for the model~\eqref{eq:CBP} with jumps as the quantity $\tilde{g}_k$ for the model~\eqref{eq:CLP} without jumps.

\smallskip

Now comes the main theoretical result of this article. The proof is in Section 3. 

\begin{thm}
\label{thm:first}
Under condition~\eqref{eq:stable-CLP}, and an additional technical condition:
\begin{equation}
\label{eq:finite-moment}
\int_{\mathbb R}z^2\,\Lambda_k(\md z) < \infty\quad \mbox{for}\quad  k = 1, \ldots, N,
\end{equation}
equivalent to $\mathbb E\left[V_j^2\right] < \infty$ for $V_j$ in~\eqref{eq:explicit-expression}, there exists a unique stationary distribution $\pi$ for the gap process $Z$ from~\eqref{eq:gap}, and the stability result~\eqref{eq:ergodic} holds. The same is true for the process~\eqref{eq:ranked-market} of ranked market weights, with another stationary distribution $\mathcal R$ instead of $\Pi$. 
\end{thm}

\begin{rmk} This statement can be extended to the case of correlated driving Brownian motions $W_1, \ldots, W_N$, and dependent L\'evy processes $L_1, \ldots, L_N$. If the $N$-dimensional L\'evy process $L = (L_1, \ldots, L_N)$ having jump measure $\Lambda$ (that is, making jumps with intensity $\lambda = \Lambda(\mathbb R^N_+)$, each jump's displacement is distributed as the normalized measure $\lambda^{-1}\Lambda(\cdot)$. Indeed, our setting in Theorem~\ref{thm:first} fits in this extended framework: If the jumps of the $k$th ranked particle are governed by a finite Borel measure $\nu_k$, $k = 1, \ldots, N$, and all jumps are independent, then
\begin{equation}
\label{eq:indep-jumps}
\Lambda = \SL_{k=1}^N\Lambda_k,\ \ \Lambda_k := \de_{0}\otimes\de_{0}\otimes\ldots\de_{0}\otimes\nu_k\otimes\de_{0}\otimes\ldots\otimes\de_{0}.
\end{equation}
\end{rmk}

\subsection{Other results} Centering $X(t)$ using~\eqref{eq:centering-procedure} produces the centered process $\overline{X} = (\overline{X}(t),\, t \ge 0)$. This is a Markov process on the hyperplane $\mathcal H := \{(x_1, \ldots, x_N) \in \BR^N\mid x_1 + \ldots + x_N = 0\}$. It is easier to prove the stability of the centered system $\overline{X}(t)$ than that of the gap process. We state this result as Theorem~\ref{thm:main}, proved in subsection 3.2, which we then use to deduce Theorem~\ref{thm:first} in subsection 3.1. Corollary~\ref{cor:LLN} is the law of large numbers, proved in subsection 3.3. For systems without jumps, it was proved in \cite{5people}.  

\begin{thm} Under conditions~\eqref{eq:stable-CLP} and~\eqref{eq:finite-moment}, the process $\ol{X}$ has a unique stationary distribution $\pi$, and the stability result~\eqref{eq:ergodic} holds. 
\label{thm:main}
\end{thm}

\begin{cor}
\label{cor:LLN}
For any bounded measurable function $f : \mathcal H \to \BR$, 
\begin{equation}
\label{eq:LLN}
\frac1T\int_0^Tf(\ol{X}(t))\,\md t \to \int_{\mathcal H}f(z)\Pi(\md z),\ \mbox{a.s. as}\ T \to \infty.
\end{equation}
In the long run, each particle $X_i$ occupies each rank $k$ on average $(1/N)$th of the time:
\begin{equation}
\label{eq:long-term-rank}
\lim\limits_{T \to \infty}\frac1T\int_0^T1\left(X_i\ \mbox{has rank}\ k\ \mbox{at time}\ t\right)\,\md t = \frac1N\quad \mbox{for each}\quad i, k = 1, \ldots, N,\quad \mbox{a.s.}
\end{equation}
\end{cor}

\subsection{Historical review} Systems of competing L\'evy particles were introduced in \cite{S2011} for $\nu_1 = \ldots = \nu_N$. Our results extend stability results \cite[Theorem 1.2(a), Theorem 1.3(a)]{S2011}. In \cite{MyOwn12}, we studied stability of systems of two competing L\'evy particles, as well as the explicit rate of exponential convergence of the gap process to its stationary distribution. Articles \cite{AB2002, Bad, KellaWhitt, Piera2008} are devoted to stability for reflected diffusions with jumps.

For competing Brownian particles (no jumps), the stationary distribution of~\eqref{eq:gap} is not known explcitily in general case. It satisfies a  complicated integro-differential equation, see \cite{5people}. For $\sigma_k = 1$, there is an explicit product-of-exponential form. Limits as $N \to \infty$ of these systems are studied in \cite{4people, S2012}. Applications in Stochastic Portfolio Theory are in \cite{JR2013b, MyOwn4}.

We know long-term behavior and scaling limits for {\it infinite systems} of competing Brownian particles in case $\sigma_k = 1$ for all $k$; see \cite{4people, DemboTsai, MyOwn6, MyOwn13, Tsai}. For the {\it infinite Atlas model:} $g_1 = 1,\, g_2 = g_3 = \ldots = 0$, there are, in fact, infinitely many (a continuum of) stationary distributions for the gap process as in~\eqref{eq:gap}, and in one of them, the bottom-ranked particle $X_{(1)}$  behaves in the long run as the fractional Brownian motion with $H = 1/4$. We can hardly hope to get such results for infinite systems with jumps, for lack of explicit formula for stationary distributions for~\eqref{eq:gap}.

\subsection{Simulations} As discussed earlier, we could not find the stationary gap distribution $\Pi$ explicitly. Instead, we perform simulations of the {\it jump Atlas model}, which is a particular case of a system of competing L\'evy particles. In this system, all particles move as Brownian motions, with $g_k = 0,\, \sigma_k = 1$, $k = 1, \ldots, N$. In addition, the (currently) bottom-ranked particle  jumps up by $b$ with intensity $\lambda$; that is, $\Lambda_1 = \lambda\delta_{b}$, and $\Lambda_k = 0$ for $k = 2, \ldots, N$.  This system has a stationary gap distribution by Theorem~\ref{thm:first}. A detailed description of simulations is given in Section 2, with code on \texttt{GitHub}. We answer the following questions:  

\begin{enumerate}[(A)]
\item In the stationary distribution, are gaps still exponential?
\item In the stationary distribution, are different gaps independent? 
\item How do mean and variance of the gap depend on jump intensity?
\item How do mean and variance depend on jump size?
\item Is the capital distribution curve still linear on the double logarithmic plot?
\end{enumerate}

The answers are: (A) No; (B) No; (C) Linear dependence (on the logarithmic scale); (D) Nonlinear dependence (on the logarithmic scale); (E) Yes.

\subsection{Acknowledgements} We are grateful to the referee for help. The second author was partially supported by NSF grants DMS 1409434 and DMS 1405210. He is grateful to \textsc{Ricardo Fernholz}, \textsc{Soumik Pal}, and \textsc{Mykhaylo Shkolnikov} for advice and discussion.

\section{Details of Simulations}

\subsection{Jump Atlas model} We simulate $300$ times a system of $N = 200$ particles with time steps $0.01$, up to time $T = 100$. We performed the following simulations: 

\begin{enumerate}
\item $\lambda = 5$ and $b \in \{0.1, 0.2, \ldots, 4.9, 5.0\}$;
\item $b = 5$ and $\lambda \in \{0.1, 0.2, \ldots, 4.9, 5.0\}$.
\end{enumerate}


\begin{figure}[t]
\centering
\subfloat[First Gap]{\includegraphics[width = 8cm, height = 6cm]{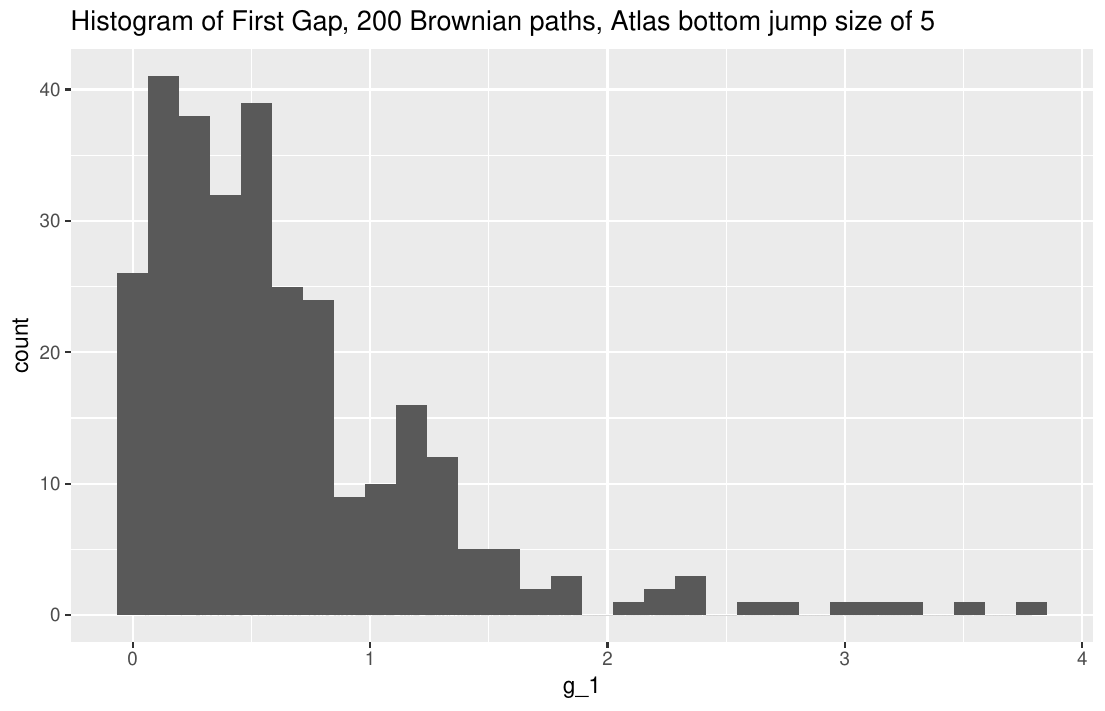}}
\subfloat[Second Gap]{\includegraphics[width = 8cm, height = 6cm]{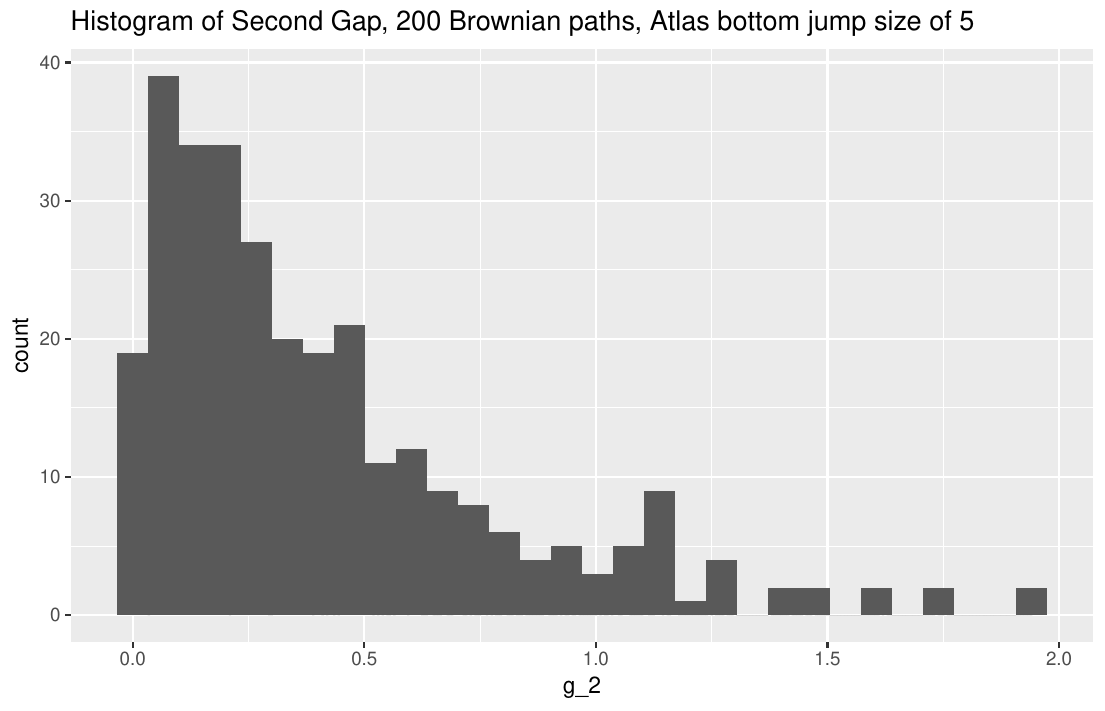}}
\caption{The first two gaps for jump size $b = 5$ and jump intensity $\lambda = 1$.}
\label{fig:hist2}
\end{figure}

\begin{figure}[t]
\center
\includegraphics[scale = 1]{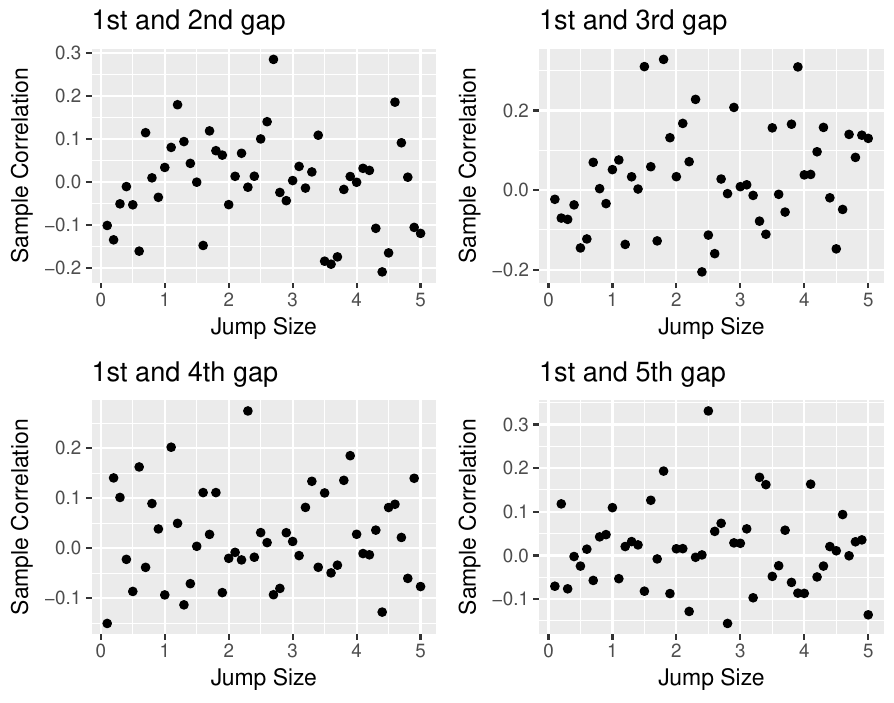}
\caption{Correlation of the first gap and the next four gaps, for jump size $b \in \{0.1, 0.2, \ldots, 4.9, 5.0\}$ and jump intensity $\lambda = 1$.}
\label{fig:corrs}
\end{figure}

\begin{figure}[t]
\centering
\subfloat[Mean vs $\lambda$, $b = 1$]{\includegraphics[width = 8cm]{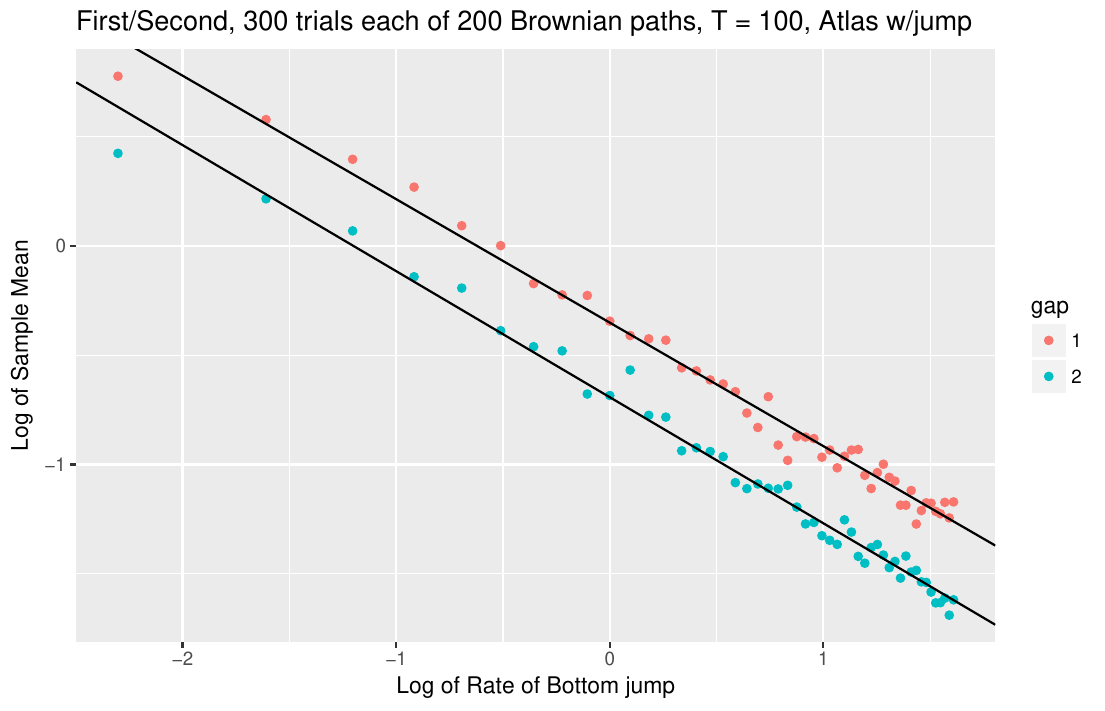}}
\subfloat[Var vs $\lambda$, $b = 1$]{\includegraphics[width = 8cm]{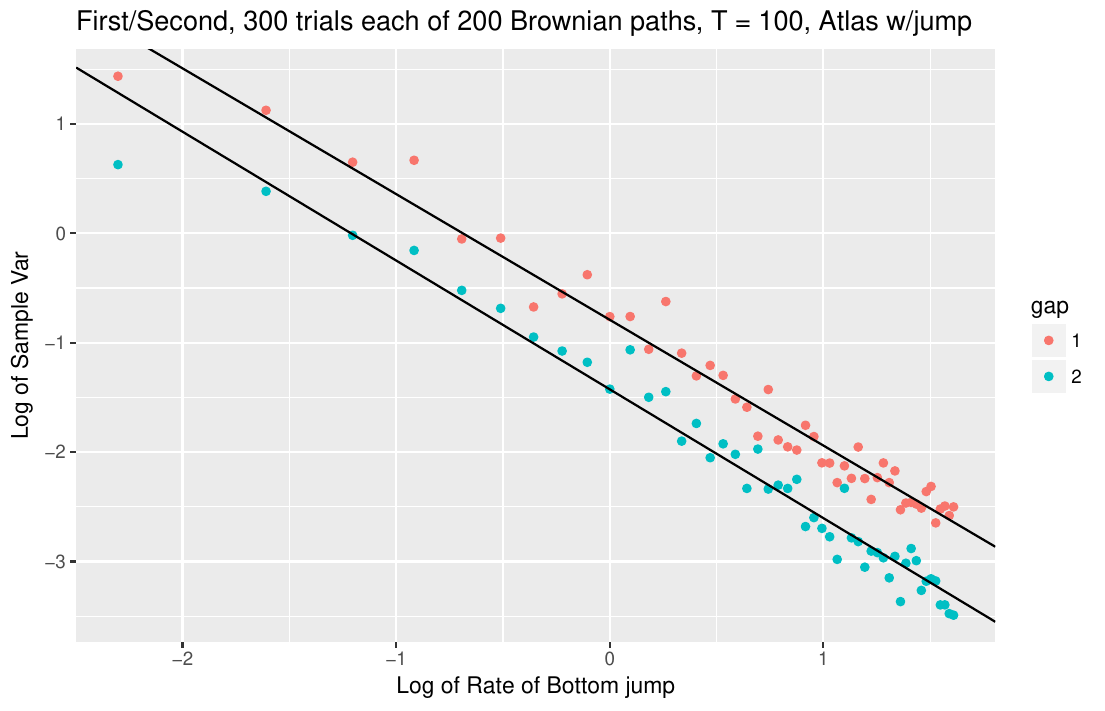}}
\\
\subfloat[Mean vs $b$, $\lambda = 1$]{\includegraphics[width = 8cm]{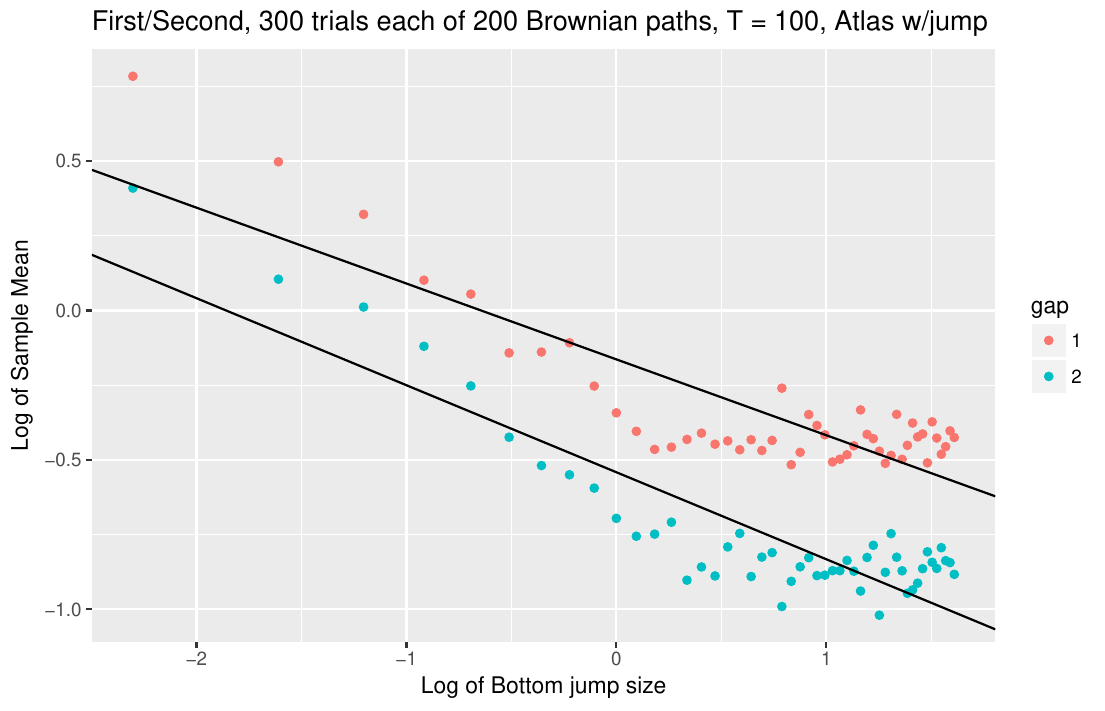}}
\subfloat[Var vs $b$, $\lambda = 1$]{\includegraphics[width = 8cm]{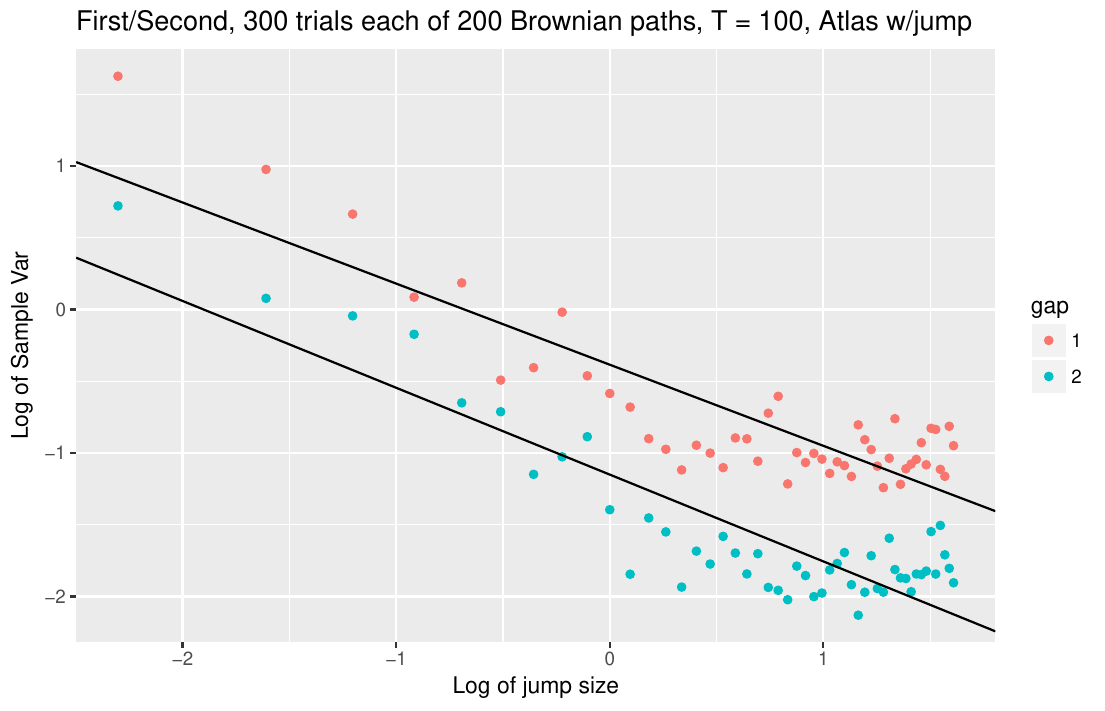}}
\caption{Log-log plot for first (red) and second (blue) gap}
\label{fig:mean-var}
\end{figure}

\begin{figure}
\subfloat[$N = 1000$, $b = 1$, $\lambda = 0.5$, $g = 0$]{\includegraphics[width = 6cm, height = 5.5cm]{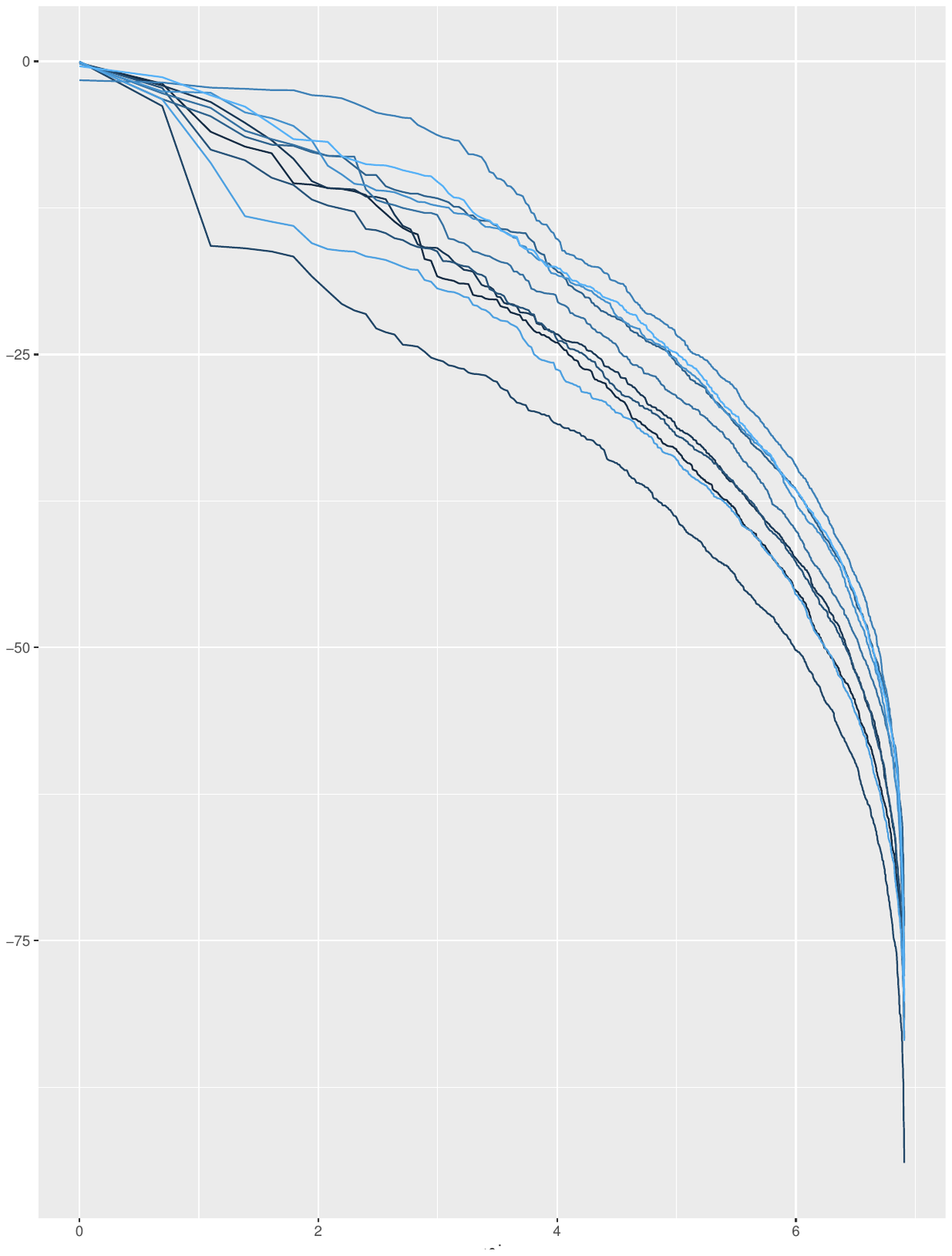}}
\subfloat[$N = 1000$, $b = 1$, $\lambda = 0.5$, $g = 5$]{\includegraphics[width = 6cm, height = 5.5cm]{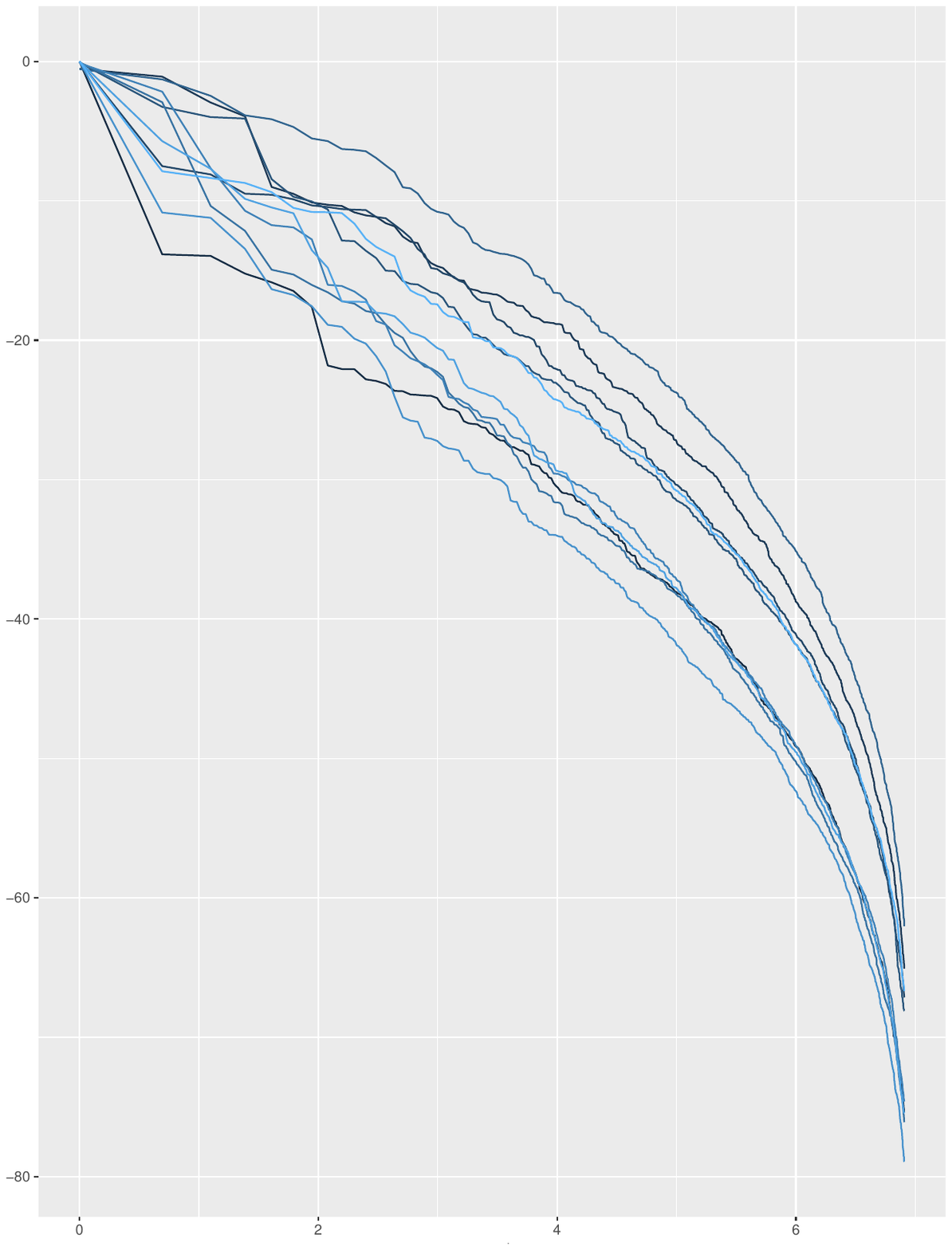}}
\subfloat[$N = 500$, $b \sim \Exp(1)$, $\lambda = 1$, $g = 0$]{\includegraphics[width = 6cm, height = 5.5cm]{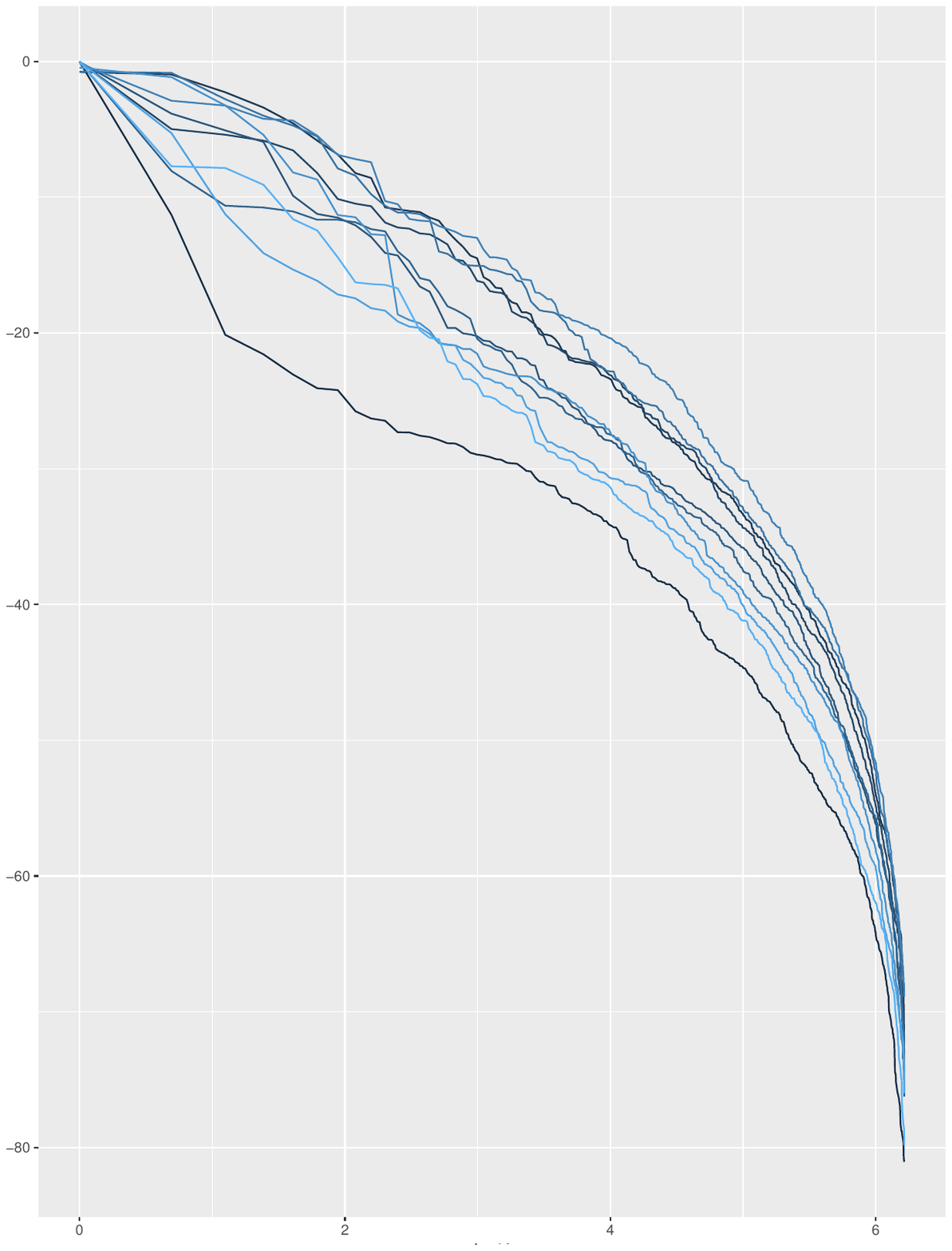}}
\caption{10 simulations of the capital distribution curve for the jump Atlas model, time steps 0.001, time horizon $T = 200$.}
\label{fig:curve}
\end{figure}

For all these 100 cases, the hypothesis that the empirical distribution is exponential, was (very much) rejected by the Kolmorogov-Smirnov test. This is in contrast with the classic Atlas model without jumps, where the stationary distribution of each gap is exponential.

\subsection{The first two gaps} We computed empirical means $\hat{\mu}_1,\, \hat{\mu}_2$ and variances $\hat{s}_1^2,\ \hat{s}_2^2$ of the first gap $Z_1 = X_{(2)} - X_{(1)}$ and the second gap $Z_2 = X_{(3)} - X_{(2)}$, and plotted them versus $\lambda \in \{0.1, 0.2, \ldots, 4.9, 5.0\}$ and $b = 1$ in Figure~\ref{fig:mean-var} \textsc{(A)} and \textsc{(B)} in a double logarithmic plot. All plots were close to linear, with almost the same slopes for the two means, and similarly for the first two variances. We did the same for $b \in \{0.1, 0.2, \ldots, 4.9, 5.0\}$ and $\lambda = 1$, and plotted empirical means and variances versus $b$ in a double logarithmic plot in Figure~\ref{fig:mean-var} \textsc{(C)} and \textsc{(D)}. These plots were not linear, but resembled functions of the type $f(x) = \max(c_1, c_2 - x)$. 

\subsection{Comparison with classic Atlas model}
Without jumps, similar linearity holds for the Atlas model from~\eqref{eq:Atlas}, if you replace $\lambda$ by $b$ from~\eqref{eq:Atlas}, \cite{5people}. The $k$th gap $Z_k \sim \Exp(r)$ with $r = 2g(1 - kN^{-1})$ has mean $m = r^{-1}$ and variance $v = r^{-2}$. Therefore, 
$$
\ln m = -\ln (g) - \ln(2) - \ln(1 - kN^{-1}),\quad \ln v = -2\ln (g) - 2\ln(2) - 2\ln(1 - kN^{-1}).
$$
That is, the double logarithmic plot of mean vs $g$ has a linear shape with slope $-1$ and intercept $-\ln(2(1 - kN^{-1}))$. For the variance vs $g$, this double logarithmic plot has a similar linear shape, but with slope $-2$ and intercept $-2\ln(2(1 - kN^{-1}))$. In this respect, the intensity 
$\lambda$  in this jump Atlas model behaves similarly to the drift $g$ in the classic Atlas model, but the jump size $b$ behaves differently.


\subsection{Correlations} We computed Pearson correlation coefficients between the first gap and the next 4 gaps. They are plotted in Figure~\ref{fig:corrs} vs jump size $b$. Although the variables are not normal, we still use Pearson correlation coefficient for the $t$-test, since the sample is large, and thus the $t$-variable (for a large degree of freedom) is close to normal. For the sample size 300, Pearson coefficient corresponding to $p = 0.05$ is $0.094$. We see that many points (much more than 5\% of the total quantity) are above this level, but not many are below $-0.094$. Thus we can rule out the hypothesis that the first and the second gaps are independent. Instead, there is evidence that they are positively correlated. This is in contrast with the case of no jumps, when the stationary distribution has independent components. 

\subsection{Capital distribution curve} Recall that in real-world stock markets and in the classic Atlas model, the capital distribution curve is linear in the double logarithmic scale. Now we simulate a jump Atlas model, modified slightly by allowing the bottom-ranked particle to have a continuous linear drift as well as jumps. Each graph displays $10$ simulations of the capital distribution curve, with $N$ the number of stocks, time steps $0.001$, time horizon $T = 200$. For the bottom-ranked particle, we denote its continuous drift by $g$, its jump size by $b$, and jump intensity $\lambda$, see Figure~\ref{fig:curve}. All other particles are driftless Brownian motions. All diffusion coefficients $\sigma_k = 1$. 

\section{Proofs} 

We introduce some notation: For a vector $x  = (x_1, \ldots, x_N) \in \BR^N$ and a permutation $\mathbf{q}$ on $\{1, \ldots, N\}$, we denote $x_{\mathbf{q}} = \left(\left(x_{\mathbf{q}}\right)_1, \ldots, \left(x_{\mathbf{q}}\right)_N\right)$, $\left(x_{\mathbf{q}}\right)_k := x_{\mathbf{q}(k)},\, k = 1, \ldots, N$. The dot product of $a = (a_1, \ldots, a_d)$ and $b = (b_1, \ldots, b_d)$ in $\BR^d$ is denoted by $a\cdot b = a_1b_1 + \ldots + a_db_d$. The Euclidean norm of $a \in \mathbb R^d$ is defined as $\norm{a} = \left[a\cdot a\right]^{1/2}$. We denote the transition kernel for the centered process $\overline{X}$ by $P^t(x, \cdot)$. 

\subsection{Derivation of Theorem~\ref{thm:first} from Theorem~\ref{thm:main}}
Theorem~\ref{thm:first} follows from Theorem~\ref{thm:main}, since the gap process $Z(t)$ is a function of the centered system $\overline{X}(t)$. Indeed, rank centered components: $\overline{X}_{(1)}(t)  \le \ldots \le \overline{X}_{(N)}(t)$, and from~\eqref{eq:centering-procedure} we see that $Z_k(t) = X_{(k+1)}(t) - X_{(k)}(t) = \overline{X}_{(k+1)}(t) - \overline{X}_{(k)}(t)$, $k = 1, \ldots, N-1$. Thus if we proved stability for the centered system $\overline{X}$, we automatically have the stability result for the gap process $Z$. 

\subsection{Proof of Theorem~\ref{thm:main}} The key ingredient is the function
\begin{equation}
\label{eq:V}
V(x) := \left[\norm{x}^2 + 1\right]^{1/2},
\end{equation}
which satisfies the following statement for some constants $b, k, r > 0$ and all $x \in \mathcal H$:
\begin{equation}
\label{eq:goal}
\ol{\CL}V(x) \le -k + b1_{\mathbb B(r)}(x),\ \mbox{where}\ \mathbb B(r) := \{x \in \mathcal H\mid \norm{x} \le r\}.
\end{equation}
We shall show~\eqref{eq:goal}, and then complete the proof of Theorem~\ref{thm:main}. The generator $\ol{\CL}$ consists of the continuous part $\ol{\CA}$ and the jump part $\ol{\CN}$. That is,  
\begin{align}
\label{eq:generator-centered}
\begin{split}
\ol{\CL}f(x) &:= \ol{\CA}f(x) + \ol{\CN}f(x),\quad f \in C^2(\mathcal H),
\\
\ol{\CA}f(x) &:= g(x)\cdot\nabla f(x) + \frac12\SL_{i=1}^N\SL_{j=1}^N\rho_{ij}(x)\frac{\pa^2f(x)}{\pa x_i\pa x_j},\\
\ol{\CN}f(x) &:= \int_{\mathcal H}\left[f(y) - f(x)\right]\,\mu_x(\md y),
\end{split}
\end{align}
where for $x \in \mathcal H$, we define $\rho_{ij}(x),\, G(x),\, \mu_x$ as follows:
$$
\rho_{ij}(x) := 
\begin{cases}
\sigma_{\mathbf{p}_x^{-1}(i)}^2,\, i = j;\\
0,\ i \ne j,
\end{cases}
\quad G(x) := \sum_{k=1}^Nx_{\mP_x(k)}\overline{g}_k,
$$
and $\mu_x$ is the push-forward to the measure $\Lambda$ with respect to the mapping
\begin{equation}
\label{eq:F-x}
F_x : \mathcal H \to \BR^N,\ \ F_x : w \mapsto x + \ol{w}_{\mP_x^{-1}}.
\end{equation}
Therefore, we can rewrite the integral in the right-hand side of $\ol{\CN}f(x)$ from~\eqref{eq:generator-centered} as
\begin{equation}
\label{eq:N-centered-new}
\ol{\CN}f(x) := \int_{\BR^N}\left[f(F_x(w)) - f(x)\right]\,\Lambda(\md w).
\end{equation}
Now, plug in $f = V$ from~\eqref{eq:V} into~\eqref{eq:generator-centered}. In \cite[Appendix, Proof of (2.18)]{BFK2005} (with notation slightly different than here), the expression $\ol{\CA} V$ is already calculated: It is the coefficient attached to $\md t$ in \cite[Appendix, Proof of (2.18), (A.13)]{BFK2005}. In our notation, we have:
\begin{equation}
\label{eq:A-V}
\ol{\CA} V(x) := \frac{x\cdot G(x)}{V(x)} + \frac{1}{V(x)}(1 - N^{-1})\sum\limits_{k=1}^N\si_k^2 - \frac{1}{V^3(x)}\SL_{k=1}^N\si_k^2x^2_{\mP_x(k)},
\end{equation}
We can rewrite~\eqref{eq:A-V} as
\begin{equation}
\label{eq:A-V-new}
\ol{\CA} V(x) := \frac{x\cdot G(x)}{V(x)} + \de_*(x),\ \ \lim\limits_{\substack{\norm{x}\to \infty}}|\de_*(x)| = 0.
\end{equation}
Consider the expression inside of the integral in~\eqref{eq:N-centered-new} for $f := V$ from~\eqref{eq:V}. The function $V$ is infinitely differentiable on $\BR^N$, and 
\begin{equation}
\label{eq:two-derivatives}
\nabla V(z) = \frac{z}{V(z)},\ \ \frac{\pa^2V}{\pa z_i\pa z_j} = \frac{\de_{ij}}{V(z)} - \frac{z_iz_j}{V^3(z)},\ i, j = 1, \ldots, N.
\end{equation}
Therefore, we have:
\begin{equation}
\label{eq:z-infty}
\lim\limits_{\norm{z} \to \infty}\left|\frac{\pa^2V(z)}{\pa z_i\pa z_j}\right| = 0,\ \ i, j = 1, \ldots, N;\quad C_V := \max\limits_{z, i, j}\left|\frac{\pa^2V(z)}{\pa z_i\pa z_j}\right| < \infty.
\end{equation}
We can write the following Taylor decomposition for all $z, u \in \BR^N$:
\begin{equation}
\label{eq:Taylor}
V(z + u) - V(z) = \frac{z\cdot u}{V(z)} + \theta(z, u).
\end{equation}
The error term $\theta(z, u)$ is given by  the following expression for some $\eta(z, u) \in [0, 1]$: 
\begin{equation}
\label{eq:error}
\theta(z, u) := \frac{\eta^2(z, u)}2\SL_{i=1}^N\SL_{j=1}^Nu_iu_j\frac{\pa^2V(z)}{\pa z_i\pa z_j}.
\end{equation}
From~\eqref{eq:z-infty} and~\eqref{eq:error}, we have
\begin{equation}
\label{eq:z-infty-2}
\lim\limits_{\norm{z} \to \infty}\theta(z, u)  = 0\ \mbox{for every}\ u \in \BR^N;\quad 
|\theta(z, u)| \le \frac{C_V}2\norm{u}^2.
\end{equation}
Letting $z := x$ and $u := \ol{w}_{\mP_x^{-1}}$ in~\eqref{eq:Taylor},
\begin{equation}
\label{eq:Taylor-new}
V(F_x(w)) - V(x) = \frac{x\cdot\ol{w}_{\mP_x^{-1}}}{V(x)} + \theta\left(x, \ol{w}_{\mP^{-1}_x}\right).
\end{equation}
From~\eqref{eq:z-infty-2} and the observation that $\norm{\ol{w}_{\mP_x^{-1}}}  = \norm{\ol{w}} \le \norm{w}$, we have
\begin{equation}
\label{eq:LDCT-new}
\left|\theta\left(x, \ol{w}_{\mP_x^{-1}}\right)\right| \le \frac{C_V}2\norm{w}^2.
\end{equation}
Combine~\eqref{eq:finite-moment},~\eqref{eq:z-infty-2}, ~\eqref{eq:Taylor-new},~\eqref{eq:LDCT-new}, and apply Lebesgue dominated convergence theorem to infer
\begin{equation}
\label{eq:N-centered-next}
\ol{\CN}V(x) = \frac1{V(x)}\int_{\BR^N}\left[x\cdot\ol{w}_{\mP_x^{-1}}\right]\,\Lambda(\md w) + \de(x),\quad \lim\limits_{\norm{x} \to \infty}\de(x) = 0.
\end{equation}
Rewrite the expression inside the integral in~\eqref{eq:N-centered-next} as follows: 
$$
x\cdot\ol{w}_{\mP_x^{-1}} = \sum_{i=1}^Nx_i\ol{w}_{\mP_x^{-1}(i)} = \sum_{k=1}^Nx_{\mP_x(k)}\ol{w}_k.
$$
Integrating with respect to $\Lambda(\md w)$, we obtain:
\begin{equation}
\label{eq:integral-lambda}
\int_{\BR^N}\left[x\cdot\ol{w}_{\mP_x^{-1}}\right]\,\Lambda(\md w) = \SL_{k=1}^Nx_{\mP_x(k)}\tilde{f}_k.
\end{equation}
Combining~\eqref{eq:N-centered-next},~\eqref{eq:integral-lambda}, we get:
\begin{equation}
\label{eq:N-centered-final}
\ol{\CN}V(x) = \frac1{V(x)}\SL_{k=1}^Nx_{\mP_x(k)}\tilde{f}_k + \de(x).
\end{equation}
As calculated in \cite[p.2302, (2.17)]{BFK2005}, via summation by parts,
\begin{equation}
\label{eq:G-drift}
x\cdot G(x) = -\SL_{k=1}^{N-1}\left(x_{\mP_x(k+1)} - x_{\mP_x(k)}\right)\left(\overline{g}_1 + \ldots + \overline{g}_k\right).
\end{equation}
Similarly, we can rewrite the sum in~\eqref{eq:N-centered-final} as 
\begin{equation}
\label{eq:N-centered-addition}
\SL_{k=1}^Nx_{\mP_x(k)}\ol{f}_k = -\SL_{k=1}^{N-1}\left(x_{\mP_x(k+1)} - x_{\mP_x(k)}\right)\left(\overline{f}_1 + \ldots + \overline{f}_k\right).
\end{equation}
Now we combine~\eqref{eq:A-V-new},~\eqref{eq:N-centered-next},~\eqref{eq:N-centered-final},~\eqref{eq:G-drift},~\eqref{eq:N-centered-addition} with the observation that $\overline{m}_j = \overline{g}_j + \overline{f}_j$. Letting $\de^*(x) := \de(x) + \de_*(x)$, we have:
\begin{equation}
\label{eq:final-formula}
\ol{\CL}V(x) = -\frac{1}{V(x)}\SL_{k=1}^{N-1}\left(x_{\mP_x(k+1)} - x_{\mP_x(k)}\right)\left(\overline{m}_1 + \ldots + \overline{m}_k\right) + \de^*(x). 
\end{equation}
Note that $V(x)/\norm{x} \to 1$ and $\de^*(x) \to 0$ as $\norm{x} \to \infty$ for $x \in \mathcal H$. As in \cite[p.2302]{BFK2005}, we have:
\begin{align}
\label{eq:upper-estimate}
\begin{split}
\SL_{k=1}^{N-1}&\left(x_{\mP_x(k+1)} - x_{\mP_x(k)}\right)\SL_{j=1}^{k-1}\overline{m}_j \le -K\norm{x},\\ K &:= N^{-1/2}\min\limits_{1 \le k \le N-1}\left(\overline{m}_1 + \ldots + \overline{m}_k\right) > 0.
\end{split}
\end{align}
From~\eqref{eq:final-formula} and ~\eqref{eq:upper-estimate}, we can prove~\eqref{eq:goal}, which is the main step of this proof of~\eqref{eq:ergodic}. We complete this proof as follows. Denote by 
$\mes$ the Lebesgue measure on $\mathcal H$. The centered system $\ol{X}$ forms a Feller continuous strong Markov process, because $X$ is a Feller continuous strong Markov process, by \cite[Theorem 2.4, Theorem 5.3, Example 1]{Sawyer}; and $\ol{X}$ is a continuous function of $X$. Next, $\ol{X}$ is a $T$-process in the terminology of \cite[Subsection 3.2]{MT1993a}: That is, there exists a nonzero function $H : \mathcal H \times \mathcal B \to [0, \infty)$ ($\mathcal B$ is the Borel $\sigma$-algebra over $\mathcal H$), lower semicontinuous in the first argument and a measure in the second argument; and for a probability measure $a$ on $[0, \infty)$, 
\begin{equation}
\label{eq:kernel} 
K_a(x, A) := \int_0^{\infty}P^t(x, A)\,a(\mathrm{d} t) \ge H(x, A),\quad A \in \mathcal B. 
\end{equation}
Let $Q^t(x, \cdot)$ be the transition kernel of a (centered) system of competing Brownian particles with the same drifts $g_k$ and diffusions $\sigma_k^2$, but without any jumps. Then 
\begin{equation}
\label{eq:Q-pve}
Q^t(x, A) > 0,\quad t > 0,\quad x \in \mathcal H,\quad A \in \mathcal B,\ \mes(A) > 0.
\end{equation}
Indeed, we treat $\mathcal H$ as $\mathbb R^{N-1}$, with  an orthogonal coordinate system. The drift coefficients are bounded, and the covariance matrix field is {\it uniformly elliptic:} There exists a $\delta > 0$ such that $\sigma(x)\xi\cdot\xi \ge \delta\norm{\xi}^2$, $\xi \in \mathbb R^{N-1}$, $x \in \mathcal H$. Under such conditions, the transition kernel has a Lebesgue density $p(t, x, y)$ which satisfies a parabolic PDE with maximum principle; thus this density is strictly positive on the whole space. This completes the proof of~\ref{eq:Q-pve}. 

Since $\lambda_0 = \lambda_1 + \ldots + \lambda_N$ is the jump intensity of the multidimensional process $(L_1, \ldots, L_N)$, we conclude that with probability $e^{-\la_0t}$ the system of competing L\'evy particles does not jump until time $t$, and behaves as a system of competing Brownian particles. Thus
\begin{equation}
\label{eq:domination}
P^t(x, \cdot) \ge e^{-\la_0t}Q^t(x, \cdot).
\end{equation}
It follows from~\eqref{eq:Q-pve} and~\eqref{eq:domination} that
\begin{equation}
\label{eq:positivity}
P^t(x, A) > 0,\quad t > 0,\quad x \in \mathcal H,\quad A \in \mathcal B,\quad \mes(A) > 0.
\end{equation}
The rest of the proof can be split in several small steps. The definition of a {\it petite set} $A$ is taken from \cite[Subsection 4.1]{MT1993a}: There exists a nontrivial measure $\nu$ on $\mathcal H$ and a probability measure $a$ on $[0, \infty)$ such that $K_a(x, \cdot) \ge \nu(\cdot)$ for all $x \in A$, with $K_a$ defined in~\eqref{eq:kernel}. 

\begin{enumerate}[(A)]

\item The discrete-time Markov chain $(\ol{X}(n))_{n \ge 0}$ is {\it irreducible} with respect to Lebesgue measure: For every subset $A \in \mathcal B$ with $\mes(A) > 0$, there is a positive probability that starting from any $x$, the Markov chain will eventually visit $A$. We took this definition from \cite[Section 4.2]{MeynBook}. This follows from~\eqref{eq:positivity}, take $t = 1$ or any positive integer.

\item We claim that all compact subsets of $\mathcal H$ are petite for the discrete-time Markov chain $(\ol{X}(n))_{n \ge 0}$. By \cite[Proposition 6.2.8(b)]{MeynBook}, this property is true for a discrete-time Markov chain which is Feller continuous (transition kernel maps bounded continuous functions to bounded continuous functions), if it has positivity property~\eqref{eq:positivity} for subset $A$ of positive $\psi$-measure, with support of $\psi$ having non-empty interior. This is indeed applicable here, since $\psi = \mes$ is the Lebesgue measure and therefore trivially has a non-empty interior.

\item Finally, \cite[Theorem 5.1]{MT1993b} implies that (in our notation and definitions): If all compact sets are petite, and if we find a function $V : \mathcal H \to [0, \infty)$ and a compact set $C \subseteq \mathcal H$ such that $\overline{\mathcal L} V(x) \le -c < 0$ for $x \notin C$, and $V$ is bounded on $C$, then the process is ergodic. We have found just such function $V$ in~\eqref{eq:V}.  
\end{enumerate}

Combine these statements (A), (B), (C), and finish the proof of~\eqref{eq:ergodic}.

\subsection{Proof of Corollary~\ref{cor:LLN}}

Let us show~\eqref{eq:LLN}. From~\eqref{eq:goal} by \cite[Theorem 4.2]{MT1993b} it follows that the process $\overline{X}$ is positive Harris recurrent: It has a unique (up to multiplication by a constant) invariant measure, and this measure is finite; it is equal to $\Pi$. Applying \cite[Theorem 17.1.7]{MeynBook} (stated and proved for discrete-time Markov chains, but similarly shown for continuous-time Markov processes), we get~\eqref{eq:LLN}. 

The limit~\eqref{eq:long-term-rank} follows from~\eqref{eq:LLN} in the same way as \cite[(2.22)]{BFK2005} follows from \cite[(2.19)]{BFK2005}: Indeed, the limit exists by~\eqref{eq:LLN}. But in this system of particles, their dynamics depends only on their current ranks. Apply any fixed permutation $\mathbf{q}$ to the particle indices. Then the  permuted system is still governed by the same law as the original system. Therefore, the limit in the right-hand side of~\eqref{eq:LLN} must be the same for all $N!$ permutations $\mathbf{q}$. The sum of these limits is $1$; therefore, each of these limits is equal to $1/N!$. There are $(N-1)!$ permutations on $\{1, \ldots, N\}$ which map $i$ to $k$. Therefore, the right-hand side of~\eqref{eq:long-term-rank} is the sum of $(N-1)!$ terms, each of which is equal to $1/N!$; the result is $N^{-1}$.

\section{Conclusion} 

We established long-term stability results for systems of competing Brownian particles with jumps, which we named {\it competing L\'evy particles}. However, we cannot find an explicit formula for stationary gap distribution, and the corresponding distribution for ranked market weights. Instead, we studied the properties of this stationary gap distribution via Monte Carlo simulations. Gaps are not independent and exponential. However, capital distribution curves are still approximately linear in the log-log scale, which captures the property of the real world markets. In the log-log scale, mean and variance of the first and second gaps depend linearly on the jump intensity, but not on the jump size. Such linearity holds for the classic Atlas model (with drift of the bottom particle instead of jump intensity).

\smallskip

Possible future research: 

\begin{enumerate}[(A)]
\item Find how the $k$th margin of the stationary distribution depends on $k$.
\item Simulate more general models of competing L\'evy particles, instead of the jump Atlas model, and replicate the research above in this more general setting.
\item Find or estmate exponential rate of convergence to the stationary distribution as $t \to \infty$.
\item Extend to competing L\'evy particles with infinitely many jumps on finite time interval. 
\end{enumerate}

\end{document}